# Pattern Formation in Annular Convection

William F. Langford [a,b], Dan D. Rusu [b]

[a] *The Fields Institute for Research in Mathematical Sciences, Toronto, Ontario M5T 3J1, Canada*

[b] *Department of Mathematics and Statistics, University of Guelph, Guelph, Ontario N1G 4S7, Canada*

**Abstract**

This study of spatio-temporal pattern formation in an annulus is motivated by two physical problems on vastly different scales. The first is atmospheric convection in the equatorial plane between the warm surface of the Earth and the cold tropopause, modeled by the two dimensional Boussinesq equations. The second is annular electroconvection in a thin smectic film, where experiments reveal the birth of convection-like vortices in the plane as the electric field intensity is increased. This is modeled by two dimensional Navier–Stokes equations coupled with a simplified version of Maxwell's equations. The two models share fundamental mathematical properties and satisfy the prerequisites for application of O(2)-equivariant bifurcation theory. We show this can give predictions of interesting dynamics, including stationary and spatio-temporal patterns.

## 1 Introduction

Pattern formation is observed in a wide variety of contexts, from slime moulds to fluid mechanics [1,2]. There is a growing realization that certain aspects of pattern formation phenomena are *model independent*, that is to say, independent of physical details of the particular model. Instead, the types of patterns which may be observed are largely determined by purely mathematical properties such as geometry and symmetry. As an illustration of this principle, we consider two different examples of pattern formation, which differ physically in the source of the instability (one thermal and the other electrical) and occur on vastly different spatial scales (ratio $\sim 10^9$). What the two examples have in common is their geometry and symmetry groups. Both take place in an annular domain and have the symmetries of the group O(2). This paper presents an analysis of those aspects of pattern formation in the two examples which are model independent. In the main, these results are qualitative, not quantitative.

Of course, a more detailed description of the behavior of a model is possible if one takes into account all of its physical properties. The extension of the qualitative analysis presented here to a detailed quantitative analysis is in progress and will be presented elsewhere. This type of quantitative analysis was first carried out in the context of the Taylor–Couette problem [3,4], which shares much with the present examples in terms of geometry and symmetry. Other examples where O(2) symmetry plays an essential role include buckling of a cylinder, flow through a cylindrical pipe and systems which are periodic on a line. The separation of quantitative from qualitative analysis is intended to illuminate the distinction between model independent phenomena and those which are intrinsically physical. The mathematical tools used here, known as Equivariant Bifurcation Theory, have been developed over the past twenty years by many mathematicians and physicists; good references are [1,2,4].

In brief outline, this paper presents two different examples of annular convection and extracts their common features on which a model independent analysis may proceed. A summary of equivariant bifurcation theory is presented, with emphasis on the key features of equivariant mode interactions in the case of O(2) symmetry. The predictions of this mathematical theory for pattern formation in annular convection are presented in tables and bifurcation diagrams. These possibilities include stationary and spatio-temporal patterns, hysteresis and mode-jumping, smooth symmetry-breaking transitions, rotating waves and standing waves.

## 2    Annular Convection: Two Problems

An annulus is the region between two concentric circles in a plane. We consider an annulus filled with an incompressible fluid, subjected to an applied force vector directed radially with magnitude independent of the angle, see Figure 1(a). If this radial force is sufficiently intense, the fluid may begin to move in interesting patterns as in Figure 1(b) and (c). A fundamental problem is how to predict these patterns.

An objection which may be raised on physical grounds is that this annulus is two dimensional, and physically no fluid is two dimensional. However, in both of the examples presented here, there are strong arguments for a two dimensional model of a three dimensional fluid. Special features of both examples, in the one case geometrical and in the other case physical, cause the fluid to behave as if it were two dimensional. This greatly reduces the com-putational difficulties in the models, and facilitates the analysis of pattern formation. Also, they provide a setting in which predictions from a two dimensional mathematical fluid model can be compared with physical reality.

## 2.1 Thermal Convection in the Equatorial Plane

Imagine that the Earth's atmosphere is spherical and ignore imperfection effects, such as continents and the seasonal tilt, which break the symmetry of a sphere rotating about a North-South axis. Then the plane slicing the Earth at the equator is a plane of symmetry and there are no forces acting on the atmosphere normal to this plane, only tangential forces. This idealized symmetry implies the possibility of flow patterns for which the equatorial plane is invariant; that is, there is no flow across the equatorial plane. For the actual Earth, it is known that this is very nearly the case.

All clouds and other weather phenomena of the atmosphere are to be found in the troposphere, which is bounded below by the warm surface of the Earth and bounded above by the cold tropopause. The equatorial plane intersects the troposphere in an annular domain. We restrict our attention to the atmospheric flow in this annulus around the equator. Recent observational work in geophysical fluid dynamics has revealed unexplained large-scale waves in the equatorial troposphere [5,6].

The geophysical fluid dynamics (GFD) in this atmospheric annulus may be modelled accurately by the Boussinesq system of equations [7]. On a thermal diffusion time scale, the system is written in the nondimensionalized form (after a translation $\tilde{T} = T - T_e$, $\tilde{p} = p - p_e$, dropping the ˜ and neglecting the Coriolis and centrifugal forces due to the rotation of the Earth)

$$\frac{\partial \vec{v}}{\partial t} + (\vec{v} \cdot \nabla) \vec{v} = \Pr \left[ -\nabla p + \nabla^2 \vec{v} + \text{Ra} \, (T \vec{\text{g}}) \right], \qquad (1)$$
$$\nabla \cdot \vec{v} = 0,$$
$$\frac{\partial T}{\partial t} + (\vec{v} \cdot \nabla) T = \nabla^2 T - T'_e(r) \, v_r,$$

where the variables are

$t$ time,

$(r, \theta)$ polar coordinates,

$\nabla$ nabla operator in polar coordinates: $\nabla \equiv (\partial/\partial r, \, r^{-1} \partial/\partial \theta)$,

$\vec{v}$ fluid velocity: $\vec{v} \equiv (v_r(r, \theta, t), v_\theta(r, \theta, t))$,

$T$ temperature variation: $T = T(r, \theta, t)$,

$T_e$ temperature of equilibrium: $T_e = T_e(r)$,

$\vec{g}$ gravity unit vector in the radial direction,

Pr Prandtl number,

Ra thermal Rayleigh number.

The boundary conditions on the two circles are

$$\vec{v} = \vec{0}, \qquad T = 0.$$

## 2.2 Electroconvection in a Smectic Liquid Crystal Film

There exist thin liquid crystal films (of type smectic-$A$) which behave like an isotropic incompressible liquid in a plane but like a plastic crystal in the normal direction [8]. The film is formed of a small integer number of planar molecular layers, of uniform thickness and behavior throughout. This structure fixes the plane of the smectic layers and suppresses any flow normal to the film. The film also is a weakly conducting dielectric, exhibiting an ohmic behavior for low electric field intensities in the plane. Experiments reported by S. Morris [9] on a thin smectic-$A$ liquid crystal film suspended in a small annulus (order of $10^{-2}$ m), with an electric field applied in the radial direction, reveal the birth of convection-like adjacent pairs of counter-rotating vortices in the plane as the electric field intensity is increased, see Figure 1.

In [10] a *surface charge model* is developed for these electrohydrodynamic (EHD) experiments, the instability being interpreted in terms of the interaction between the exterior three dimensional electric field (just outside the film) and the surface charges developed on the smectic-$A$ film viewed as a suspended surface with two (top and bottom) faces. In this model, it is the discontinuity in the normal component of the electric field across the surfaces of the smectic film which creates surface charges. These surface charges in turn experience a force due to the tangential component of the electric field, which drives the convection instability.

On the other hand, a *volume charge model* may be constructed, taking into account the spontaneous generation of ions in the bulk of the liquid crystal through electrochemical dissociation/recombination, perhaps involving dopant, (see [11])

$$AB \rightleftharpoons A^- + B^+.$$
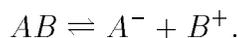

This dissociation by itself yields equal numbers of positive and negative ions and a zero net charge; however, near the positive and negative electrodes respectively positive and negative net charges build up, due to injection and/or annihilation mechanisms. Locally (near the electrodes) diffusion is important in this model. A convective instability arises from the interaction between these net volume charges and the two-dimensional electric field tangent to the film, analogous to that in the Boussinesq model. Experiments reported in [8] on rectangular smectic-A films, with and without electrode contact, suggest that some contact mechanism at the electrodes, perhaps electrochemical, is required for pattern formation. The following equations represent a volume charge model which neglects the surface charge mechanism, and strongly resembles the above Boussinesq approximation in GFD. After scaling on a conduction charge decay time, it may be written in the nondimensionalized form (after a translation $\tilde{q} = q - q_\mathrm{e}$, $\tilde{\vec{E}} = \vec{E} - \vec{E}_\mathrm{e}$, $\tilde{p} = p - p_\mathrm{e}$ and dropping the $\tilde{\ }$) as

$$\frac{\partial \vec{v}}{\partial t} + (\vec{v} \cdot \nabla)\vec{v} = \Pr\left[-\nabla p + \nabla^2 \vec{v} + \mathrm{Ra}\left(q\vec{E} + q_\mathrm{e}\vec{E} + q\vec{E}_\mathrm{e}\right)\right], \qquad (2)$$

$$\nabla \cdot \vec{v} = 0,$$

$$\frac{\partial q}{\partial t} + (\vec{v} \cdot \nabla)q = q + \delta \Pr \nabla^2 q - q'_\mathrm{e}(r)\, v_\mathrm{r},$$

$$q = \nabla \vec{E},$$

$$\vec{E} = -\nabla \phi,$$

where the variables are

$q$ net charge density variation: $\quad q = q(r, \theta, t)$,

$q_\mathrm{e}$ net charge density at equilibrium: $\quad q_\mathrm{e} = q_\mathrm{e}(r)$,

$\vec{E}$ in-plane electric field variation: $\quad \vec{E} \equiv (E_\mathrm{r}(r, \theta, t), E_\theta(r, \theta, t))$,

$\vec{E}_\mathrm{e}$ in-plane electric field of equilibrium: $\quad \vec{E}_\mathrm{e} \equiv (E_\mathrm{e}(r), 0)$,

$\phi$ electric potential,

Pr electric Prandtl number,

Ra electric Rayleigh number,

$\delta$ dimensionless ratio of charge diffusion coefficient to kinematic viscosity.

The boundary conditions on the two circles are

$$\vec{v} = \vec{0}, \qquad q = \phi = 0.$$

If two or more species of ions with different mobilities are important, then this model may be expanded with a different charge continuity equation for each species of ion. In the forcing term of the Navier–Stokes equation the overall total net charge will appear.

*2.3  Common Features of the two Models*

In spite of their obvious differences, the two annular convection models from GFD (geophysical fluid dynamics) and EHD (electrohydrodynamics) have many common features. Both are especially interesting because they provide physically relevant examples of two dimensional Navier–Stokes flows. Both have O(2)-symmetry, i.e. the symmetry of the circle under rotations and reflections across any diameter.

There is a destabilizing radial force, due to gravity and the temperature-induced density variation for the GFD model and due to the in-plane electric field intensity and the electric charge density variation for the EHD model. Each of these destabilizing forces is opposed by a a stabilizing force due to viscosity.

In both cases there are trivial (pre-convection) solutions which are known exactly (analytically). The neutral stability boundary is fairly easy to calculate (numerically) after an exact linearization about the known trivial solution. The initial loss of stability (bifurcation) is to a convection mode, which can be described as a chain of pairs of clockwise/counterclockwise vortices strung periodically around the annulus, see Figure 1. The number of such pairs is the mode number $m$, and this depends on the radius ratio $\eta \in (0, 1)$. The two important parameters are the Rayleigh number Ra (depending on temperature or voltage respectively, as appropriate to the model) and the radius ratio $\eta$. Codimension-two points for two different mode numbers exist in the parameter space. We show that the prerequisites for the application of O(2)-Equivariant Bifurcation Theory are satisfied for both models. This implies the possibility of rich mode interactions.

## 3  Equivariant Bifurcation Theory

We investigate these two problems in parallel, since the mathematical tools [1–4] for both are essentially the same. Consider a parametrized family of differential equations, representing for example (1) or (2),

$$y' = f(y, \mu), \quad f : \mathcal{X} \times \mathcal{R}^p \to \mathcal{Y}, \tag{3}$$

where $f$ is a smooth $C^k$ function ($k > 2$) such that $f(0,0) = 0$, $\mathcal{X}, \mathcal{Y}$ are appropriate function spaces and $\mathcal{R}^p$ is the parameter space ($\mathcal{R}$ = set of real numbers). Local bifurcation theory deals with the question: What qualitative changes can happen to the set of solutions near $y = 0$, as the parameter $\mu$ varies near zero? The answer depends in large part on the eigenvalues of the linearization of $f$ at the origin, i.e. $D_y f(0, \mu)$. Typically there are three possibilities: either there is no qualitative change (the eigenvalues remain hyperbolic as $\mu$ varies), or there is a steady-state bifurcation (a real eigenvalue satisfies $\lambda(\mu) = 0$ for some $\mu$), or there is a Hopf bifurcation of a limit cycle (a complex eigenvalue satisfies $\text{Re}\,[\lambda(\mu)] = 0$ for some $\mu$). In this paper, we consider only the case of steady-state bifurcation, but we allow the possibility of multiple real eigenvalues vanishing for the same value of $\mu$, i.e. of mode-interactions. Without loss of generality, we translate $\mu$ so that this multiple zero eigenvalue occurs at $\mu = 0$, i.e. $\lambda(0) = 0$.

### 3.1 Reduction to a Center Manifold

Suppose that this zero eigenvalue $\lambda(0) = 0$ has finite multiplicity and that its complement in the spectrum of $D_y f(0,0)$ is bounded away from the imaginary axis, in the negative half-plane. Then, under reasonable conditions on $f$, which are satisfied for many PDE of the general form (3), there exists a center manifold in $\mathcal{X}$ which is tangent to, and of the same dimension as, the (generalized) eigenspace of this zero eigenvalue [12]. Even in cases where the existence of this center manifold has not yet been proven rigorously, it is possible to calculate it in a formal manner. In this paper, we set aside the pure existence question and assume there is a normally hyperbolic center manifold for the annular convection models above. Then, locally, all solutions off the center manifold approach it exponentially with increasing time. Thus the study of the long-time behavior of solutions of (3) near the origin is reduced to that of solutions on this finite-dimensional manifold. This reduction to finite dimensions has been calculated, for example, in the case of the Taylor–Couette problem in [3,4].

Furthermore, this reduction process preserves the symmetry properties of the system, as described in the next subsection.

### 3.2 Symmetry: Definitions and Tools

Consider a system of ordinary differential equations

$$y' = f(y, \mu), \quad y \in \mathcal{R}^n, \quad \mu \in \mathcal{R}^p, \tag{4}$$

and let $\Gamma$ be a compact group acting linearly on $\mathcal{R}^n$ i.e.

$$\Gamma: \mathcal{R}^n \to \mathcal{R}^n, \quad \Gamma \subset \mathrm{GL}(n),$$

where $\mathrm{GL}(n)$ is the group of all invertible $n \times n$ matrices over $\mathcal{R}$. Then, system (4) is said to have the symmetry $\Gamma$ (or is $\Gamma$-*equivariant*) if

$$f(\gamma y, \mu) = \gamma f(y, \mu), \quad \forall\, \gamma \in \Gamma, \ (y, \mu) \in \mathcal{R}^n \times \mathcal{R}^p. \tag{5}$$

In this paper we are concerned with the *orthogonal group* $\Gamma = \mathrm{O}(2)$ on $\mathcal{R}^2$ which is generated by $\mathrm{R}_\theta$, the planar rotations by $\theta$, and the flip $\kappa$, the reflection across the $x$-axis. Then, the definition (5) says that $f$ has symmetry $\mathrm{O}(2)$ if $f$ commutes with every $\gamma \in \mathrm{O}(2)$. This algebraic identity agrees with our geometrical notion of symmetry. Indeed, taking $\gamma = \mathrm{R}_\theta$, the definition of $\Gamma$-equivariance says that the value of $f$ at the point obtained by rotating $y$ through angle $\theta$ is the same as $f(y)$ rotated through $\theta$. Similarly, taking $\gamma = \kappa$, it says that $f$ at the reflection of $y$ is the reflection of $f(y)$.

For the remainder of this section, let us temporarily suppress the parameter $\mu$ in (4) and (5), since $\mu$ plays no direct role in the symmetry; that is rewrite (5) as

$$f(\gamma x) = \gamma f(x), \quad \forall\, \gamma \in \Gamma, \ y \in \mathcal{R}^n.$$

The *group orbit* of a point $x \in \mathcal{R}^n$, with respect to $\Gamma$, is defined as

$$\Gamma x =: \{\gamma x \mid \gamma \in \Gamma\}. \tag{6}$$

Suppose that $x$ is a steady-state, i.e. $f(x) = 0$. Then $0 = \gamma 0 = \gamma f(x) = f(\gamma x)$, for all $\gamma \in \Gamma$, so the entire group orbit is composed of steady-states. We have the following dichotomy (and both cases are interesting): either $\gamma x = x$ (i.e. $x$ has the symmetry $\Gamma$) or $\gamma x \neq x$ (i.e. $\gamma x$ is a new solution). In the latter case, each $\gamma x$ is called a *conjugate* of x. Conjugates have the same dynamic behavior. Similar considerations may be employed for limit cycles and other types of solutions of (4).

The *isotropy subgroup* of $x \in \mathcal{R}^n$ is defined as

$$\Sigma_\mathrm{x} =: \{\gamma \in \Gamma \mid \gamma x = x\} \tag{7}$$

and represents the total symmetry of $x$. The *fixed-point subspace* of $\Sigma$, defined by

$$\mathrm{Fix}(\Sigma) =: \{x \in \mathcal{R}^n \mid \gamma x = x, \ \forall\, \gamma \in \Sigma\} \tag{8}$$

represents the set of all points fixed by the subgroup $\Sigma$. It is a subspace of $\mathcal{R}^n$, which is invariant for any $\Gamma$-equivariant $f$, i.e.

$$f(\mathrm{Fix}(\Sigma)) \subset \mathrm{Fix}(\Sigma). \tag{9}$$

Even though $f$ is nonlinear, we can restrict the problem to lower dimensional linear subspaces on which $f$ is invariant. This is much like the decomposition of a linear problem into eigenspaces. Thus (9) yields a very powerful tool in the analysis of symmetric systems.

Next we show how symmetry can lead to multiple eigenvalues. Assume that 0 is a steady-state, i.e. $f(0) = 0$. Let $L =: (\partial f/\partial y)(0)$ (the linearization of $f$ at 0), and suppose that $L$ has the eigenvalue $\lambda = 0$. Then $V =: \ker L \neq \{0\}$. By differentiating the equivariance relation $f(\gamma x) = \gamma f(x)$, we obtain, at $x = 0$, that $L\gamma = \gamma L$, so $L$ and $\gamma$ are commuting linear operators. Thus $\Gamma$ leaves $V$ invariant. For most $\Gamma$, we must have $\dim V > 1$ for a faithful representation on $V$. For instance, in order to have the full symmetry of $\Gamma = \mathrm{O}(2)$ (the symmetry of the annulus) or $\Gamma = \mathrm{D_n}$ (the $n$-gon symmetry), we must have at least $\dim V = 2$. Then $\lambda = 0$ is a multiple eigenvalue "forced" by the symmetry.

*3.3 Normal Forms*

For $\Gamma$ compact, $\Gamma$-equivariant maps $f$ can be written in a very useful simplified form. First consider the class of *invariant functions*, i.e. smooth functions $g \colon \mathcal{R}^n \to \mathcal{R}$ satisfying

$$g(\gamma x) = g(x) \quad \forall \, \gamma \in \Gamma, \ x \in \mathcal{R}^n. \tag{10}$$

Then there exists a finite set of invariant polynomials $\{u_1, ..., u_s\}$ with $u_k(\gamma x) = u_k(x)$ for all $\gamma \in \Gamma$ and $k \in \{1, ..., s\}$, such that any invariant function $g$ can be written as $g(x) = h(u_1(x), ..., u_s(x))$ for some smooth function $h$. This result holds for $g$ and $h$ polynomials or germs. The set $\{u_1, ..., u_s\}$ is called a *Hilbert basis*. Next follows the Equivariant Normal Form Theorem: for a given equivariant $f \colon \mathcal{R}^n \to \mathcal{R}^n$ (polynomial or germ) there exists a finite set of $\Gamma$-equivariant polynomial maps $\overline{v}_1, ..., \overline{v}_r \colon \mathcal{R}^n \to \mathcal{R}^n$ such that

$$f(x) = \sum_{k=1}^{r} h_k(u_1(x), ..., u_s(x)) \, \overline{v}_k(x) \tag{11}$$

with $h_k$ and $u_j$ as above. Equation (11) gives the *equivariant normal form* of $f$. For the case of $\Gamma = \mathrm{O}(2)$, the sets of basic invariant polynomials and equivariants were obtained by Buzano and Russo, see [13,14].

*3.4 Mode Interactions*

Roughly speaking, the *codimension* of a bifurcation is the number of equations (in the parameters) which must be satisfied in order for that bifurcation to occur, under the assumption that certain nondegeneracy conditions hold. Typical bifurcations have codimension one. A simple setting in which a codimension two bifurcation can occur is as follows. In a two-parameter system, as we vary the two parameters $\mu \equiv (\mu_1, \mu_2) \in \mathcal{R}^2$ in $y' = f(y, \mu)$, we encounter codimension-one curves of neutral stability, each corresponding to a zero eigenvalue (or complex eigenvalue with zero real part) with the symmetry of the corresponding kernel (see Figure 2). The curves corresponding to bifurcations of codimension-one are well-understood. Generically these curves will intersect transversally at isolated points (codimension-two bifurcation points) where *mode interactions* occur. Near a codimension-two point as in Figure 2 (at the intersection of two codimension-one curves) one naturally expects to find both of the codimension-one behaviors, and generally does so. However, the interesting fact is that often there is much more to be found, since nonlinear "mode interactions" can produce interesting new phenomena.

## 4 Analysis of Mode Interactions

The analysis of the mode interactions takes place in the following steps. For more details, see [1].

• Adjust the parameters to a neighborhood of the codimension-two point of O(2) mode interactions, i.e where two eigenvalues are simultaneously zero. This involves only linear stability analysis.

• Perform the *center manifold reduction* from the original PDE to a four dimensional ODE system

$$y' = g(y, \mu), \quad g : \mathcal{R}^4 \times \mathcal{R}^2 \to \mathcal{R}^4, \tag{12}$$

where

$$g(0, \mu) \equiv 0, \quad \mathrm{D}g(0, 0) = 0, \quad g(\gamma y, \mu) = \gamma g(y, \mu), \quad \forall\ \gamma \in \mathrm{O}(2) \tag{13}$$

(we assume O(2) acts faithfully on both modes, which implies a four dimensional system).

- Identify the state space $\mathcal{R}^4$ with $\mathcal{C}^2$ using the complex coordinates $(z_1, z_2)$. The action of O(2) on $\mathcal{C}^2$ is given by $\theta(z_1, z_2) = (e^{im\theta} z_1, e^{il\theta} z_2)$ (rotation) and $\kappa(z_1, z_2) = (\overline{z}_1, \overline{z}_2)$ (reflection) where $l$ and $m$ are the integer mode numbers. Based on available experimental evidence, we may assume that $l = m + 1$.

- The Hilbert basis for the ring of O(2)-invariants as obtained in [13] is $\{u, v, w\}$ given by

$$u =: |z_1|^2, \quad v =: |z_2|^2, \quad w =: z_1^{m+1} \overline{z}_2^m + \overline{z}_1^{m+1} z_2^m. \tag{14}$$

- The generators of the module of O(2)-equivariants (over the above ring of invariants) from [13] is

$$\left\{ \begin{bmatrix} z_1 \\ 0 \end{bmatrix}, \begin{bmatrix} 0 \\ z_2 \end{bmatrix}, \begin{bmatrix} \overline{z}_1^m z_2^m \\ 0 \end{bmatrix}, \begin{bmatrix} 0 \\ z_1^{m+1} \overline{z}_2^{m-1} \end{bmatrix} \right\}. \tag{15}$$

- It follows that the normal form for the O(2)-equivariant system is

$$\frac{d}{dt} \begin{bmatrix} z_1 \\ z_2 \end{bmatrix} = p_1 \begin{bmatrix} z_1 \\ 0 \end{bmatrix} + p_2 \begin{bmatrix} 0 \\ z_2 \end{bmatrix} + q_1 \begin{bmatrix} \overline{z}_1^m z_2^m \\ 0 \end{bmatrix} + q_2 \begin{bmatrix} 0 \\ z_1^{m+1} \overline{z}_2^{m-1} \end{bmatrix} \tag{16}$$

together with the two equations for $\overline{z}_1, \overline{z}_2$ which are complex conjugates of these. Here $p_1, p_2, q_1, q_2$ are real functions of $(u, v, w, \mu)$, and satisfy

$$p_1(0) = p_2(0) = 0.$$

- Compute the isotropy subgroups of O(2). They are conveniently organized in the form of a *lattice of isotropy subgroups*, where the arrows represent inclusions:

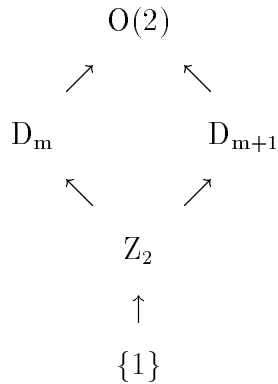

Table 1: Lattice of Isotropy Subgroups.

- In Table 1, $Z_2$ is the reflection symmetry and $\{1\}$ has no symmetry. Each isotropy subgroup determines a corresponding fixed-point subspace.

- Restrict the system to each fixed-point subspace and solve it. We recall that each fixed-point subspace is invariant for the flow (i.e. for time-evolution of solutions) and the system simplifies when restricted to each fixed-point subspace, as follows:

| Isotropy Subgroup $\Sigma$ | Fixed-Point Subspace Fix($\Sigma$) | Equations Restricted to Fixed-Point Subspace $g|_{\text{Fix}(\Sigma)} = 0$ |
|---|---|---|
| O(2) | $\{0\}$ | $0 = 0$ |
| $D_m$ | $\{(x, 0)\,|\,x \in \mathcal{R}\}$ | $x\,p_1(x^2, 0, 0, \mu) = 0$ <br> $z_1 = x \in \mathcal{R}$ <br> $z_2 = 0$ |
| $D_{m+1}$ | $\{(0, y)\,|\,y \in \mathcal{R}\}$ | $y\,p_2(0, y^2, 0, \mu) = 0$ <br> $z_1 = 0$ <br> $z_2 = y \in \mathcal{R}$ |
| $Z_2$ | $z = (x, y)$ <br> $x, y \in \mathcal{R}$ | $h_1 \equiv p_1 + q_1 x^{m-1} y^m = 0$ <br> $h_2 \equiv p_2 + q_2 x^{m+1} y^{m-2} = 0$ |
| $\{1\}$ | $\mathcal{C}^2$ | $p_1 = p_2 = 0$ <br> $q_1 = q_2 = 0$ |

Table 2: Restrictions to Fixed-Point Subspaces.

- Compute numerical values of the leading coefficients of (low degree) nonlinear terms in $p_1, p_2, q_1, q_2$ and thus determine local dynamics in each of the fixed-point subspaces in this table. (This is work in progress.)

## 5  Bifurcation Scenarios

A general analysis of the O(2)-equivariant normal form (16) has been performed by Dangelmayr [14] and reviewed in [1]. A rich assortment of behaviors is possible. The determination of which behaviors actually occur in a given physical system depends on the specific values of certain coefficients (to be computed). Here we describe some interesting possibilities which follow from the general analysis, assuming $m \geq 2$. (The case $m = 1$ is quite different and will not be discussed here.)

First, consider the equations for steady-states in the subspace Fix($Z_2$), i.e. $g|_{\text{Fix}(Z_2)} = 0$:

$$p_1 x + q_1 x^m y^m = 0, \qquad (17)$$
$$p_2 y + q_2 x^{m+1} y^{m-1} = 0,$$

where $p_1, q_1, p_2, q_2$ are functions of $u = x^2, v = y^2, w = x^{m+1} y^m$ and the parameters. It can be shown that, generically, there are no steady-states other than those in Fix($Z_2$). Thus symmetry has reduced the dimension of the problem from four to two. It is clear that (17) has the trivial solution $(x, y) = (0, 0)$. We seek nontrivial solutions.

Suppose now that the mode number $m$ is very large, as it would be in the GFD example (1) where the radius ratio is very close to 1. Then locally, the two terms including $q_1, q_2$ are very small compared to the other two terms involving $p_1, p_2$ (which are independent of $m$). Similarly, among the invariants, $w = \text{Re}\left(\bar{z}_1^{m+1} z_2^m\right)$ is very small for large $m$ compared with $u$ and $v$. Thus, to obtain a first approximation when $m$ is large, we may set $w, q_1, q_2$ all equal to zero. The resulting equations are the well-known bifurcation equations for the case of $Z_2 \oplus Z_2$ symmetry, which arise in cases of pitchfork/pitchfork, pitchfork/Hopf and nonresonant Hopf/Hopf mode interactions (with no additional symmetry or degeneracy). The bifurcation equations for this case to leading order are

$$x\left(\alpha + p_{11} x^2 + p_{12} y^2\right) = 0, \qquad (18)$$
$$y\left(\beta + p_{21} x^2 + p_{22} y^2\right) = 0,$$

where $\alpha, \beta$ are the unfolding parameters and $p_{ij}$ are the obvious derivatives, evaluated at the origin. We assume the nondegeneracy conditions $p_{11} p_{22} \neq 0$ and $p_{11} p_{22} - p_{12} p_{21} \neq 0$. The complete set of bifurcation diagrams for this situation has been published in [15], see also [16,1]. Two interesting cases are shown in Figure 3, where the bifurcation parameter $\mu$ is related to the unfolding parameters by

$$\alpha = \mu, \quad \beta = \mu + \delta,$$

and $\delta > 0$ is a small "splitting parameter". In all cases, there are two primary pitchfork branches from the trivial solution, given by

$$x^2 = -\frac{\mu}{p_{11}}, \quad y = 0, \qquad (19)$$

$$x = 0, \quad y^2 = -\frac{\mu + \delta}{p_{22}},$$

and secondary "mixed-mode" branches with both of $x, y \neq 0$, which are solutions of the linear system in $x^2, y^2$ obtained from (18) on dividing out $x \neq 0$ and $y \neq 0$. The existence of all these branches persists when the higher order terms are restored to the equations (18), by the Implicit Function Theorem. In the case of $Z_2 \oplus Z_2$ symmetry (with $w, q_1, q_2$ suppressed) the symmetries of these branches also are preserved: the primary branches remain symmetric pitchforks, and the secondary branches occur as four conjugate $(\pm x, \pm y)$ under the $Z_2 \oplus Z_2$ symmetry. It is convenient in this case to identify these conjugates and draw them as single branches in a planar bifurcation diagram with amplitude plotted against bifurcation parameter, see Figure 3.

Now let us return to the equations for steady-states of the original O(2)-symmetric system given by (17). The $m$-dependent terms break the $Z_2 \oplus Z_2$ symmetry but the existence of all the branches in the $Z_2 \oplus Z_2$ symmetric case is preserved. For sufficiently large $m$, the symmetry-breaking effects will be very small and the bifurcation diagram will be indistinguishable from the $Z_2 \oplus Z_2$ case. For all values of $m \geq 2$, we note that (17) retains a $Z_2$ equivariance which depends on the parity of $m$:

$$(x, y) \to (-x, y), \quad m \text{ odd}$$
$$(x, y) \to (x, -y), \quad m \text{ even}.$$

For definiteness, let us assume $m$ is odd (the even case is similar). Then in the bifurcation diagrams of Figure 3 (full 3D versions), the $Z_2$-symmetry: $(x, y) \to (-x, y)$ implies that the $D_m$-mode primary pitchfork is preserved and the secondary bifurcations occur at the same points on the two prongs of this pitchfork. However, the lack of a $Z_2$-symmetry: $(x, y) \to (x, -y)$ means that the four secondary branches occur in two distinct sets of conjugate pairs. Thus in drawing the conjugacy classes of branches in this case we show two distinct $Z_2$-symmetric (mixed-mode) branches, as in Figure 4. A further implication of the breaking of the $Z_2$-symmetry: $(x, y) \to (x, -y)$ is that the secondary bifurcation from the $D_m$-symmetric $x$-mode becomes an imperfect pitchfork, which implies the possibility of a transcritical bifurcation with a nearby saddlenode (or limit point) bifurcation. This is indicated in Figure 4. There is an additional "hidden symmetry" which preserves only the primary pitchfork for the $D_{m+1}$-symmetric $y$-mode with $x = 0$, see [1] for details.

The bifurcation diagrams in Figure 3 show secondary (mixed-mode) branches which link the two primary (pure-mode) branches, in two different ways. In the case (a), the secondary branch is itself unstable but it effects an exchange of stabilities between the primary branches (stable solutions are indicated by heavy lines). If the parameter $\mu$ is varied up and down through this range, hysteretic jumps between the two primary branches will occur. In case (b) of Figure 3, the secondary branch is stable and it effects a smooth stable transition between the primary branches. Along the stable secondary branch joining the $D_m$ to the $D_{m+1}$ primary branch, the corresponding picture in the annulus has a small counter-rotating vortex pair created near the outer circle between two existing vortices. This new pair grows and eventually matches the other pairs as the $D_{m+1}$ branch is approached. On this secondary branch the symmetry is broken, from $D_m$ (or $D_{m+1}$) to $Z_2$.

Additional phenomena can be predicted if one of the primary bifurcations is subcritical with the other remaining supercritical, as in Figure 4. Then, for the $Z_2 \oplus Z_2$ symmetric normal form, as shown in [16,15], there is generically a Hopf bifurcation of a periodic solution from the mixed-mode branch. In fact, there are four identical Hopf bifurcations but all are conjugate under the $Z_2 \oplus Z_2$ symmetry. Returning now to equations (17) for the steady-states of our original O(2)-symmetric system, these Hopf bifurcations are preserved but now exist as two distinct pairs of $Z_2$-conjugates. Thus, under these conditions we predict the Hopf bifurcation of two distinct branches of periodic solutions. Because these are purely temporal oscillations (not spatio-temporal) we may call them *standing wave* solutions. They are marked SW in Figure 4.

Still another type of solution is possible, as a consequence of the original O(2)-symmetry. The primary (pure-mode) and secondary (mixed-mode) solution branches found for the reduced system (17) actually exist as full "circles" of solutions, conjugate under the SO(2)-symmetry, for the original equations. Dangelmayr [14] pointed out that in addition to the two eigenvalues associated with the 2D-systems above, there is one eigenvalue identically zero due to this SO(2)-symmetry, and a fourth eigenvalue which may change sign. When this happens, the result is bifurcation of a *rotating wave*, involving a slow motion along the SO(2) group orbit. Here by rotating wave is meant a spatio-temporal pattern with the property that there is a rotating coordinate frame in which the pattern becomes stationary. This rotating wave is marked RW in Figure 4.

## 6 Conclusions

A rich variety of stationary and spatio-temporal patterns is possible in systems with the symmetries of annular convection problems. In the geophysical fluid dynamics problem (1), the mode number $m$ is very large, which implies that the bifurcation diagrams are practically indistinguishable from the classical $Z_2 \oplus Z_2$ symmetric case [16,15]. Either hysteresis or smooth transitions between modes, and stationary waves are possible behaviors. In the electroconvection problem (2), mode numbers $m$ are in the single-digit range.

This allows even richer behaviors, for example the possibility of a rotating wave. Further calculations (in progress) will determine which of these can be predicted in given experiments.

## Acknowledgement

The authors would like to thank the participants at the CPiP'98 Conference for their discussions and insightful questions. Answers to these questions have been incorporated in the final version of this paper.

## References

[1] M. Golubitsky, I. Stewart and D.G. Schaeffer, Singularities and Groups in Bifurcation Theory II (Springer–Verlag, New York, 1988).

[2] M.C. Cross and P.C. Hohenberg, Rev. Mod. Phys. 65 (1993) 851.

[3] M. Golubitsky and W.F. Langford, Physica D 32 (1988) 362.

[4] P. Chossat and G. Iooss, The Taylor–Couette Problem (Springer–Verlag, New York, 1994).

[5] T. Tsuda, Y. Murayama, H. Wiryosumarto, S.W.B. Harijono and S. Kato, J. Geophysical Research 99, D5 (1994) 491 and 507.

[6] K. Hamilton, J. Geophysical Research 102, D12 (1997) 457.

[7] P.G. Drazin and W.H. Reid, Hydrodynamic Stability (Cambridge University Press, Cambridge, 1981).

[8] S.W. Morris, J.R. de Bruyn and A.D. May, J. Stat. Phys. 64 (1991) 1025.

[9] Z.A. Daya, V.B. Deyirmenjian, S.W. Morris and J.R. de Bruyn, Phys. Rev. Lett. 80 (1998) 964.

[10] Z.A. Daya, S.W. Morris and J.R. de Bruyn, Phys. Rev. E 55 (1997) 2682.

[11] S.W. Morris, J.R. de Bruyn and A.D. May, Phys. Rev. A 44 (1991) 8146.

[12] A. Vanderbauwhede and G. Iooss, Dynamics Reported 1 *New series* (1992) 125.

[13] E. Buzano and A. Russo, Annali di Matematica Pura ed Applicata (IV) 146 (1987) 217.

[14] G. Dangelmayr, Dynam. Stab. Syst. 1 (1986) 159.

[15] W.F. Langford and G. Iooss, Bifurcation Problems and their Numerical Solutions (eds H.D. Mittelmann and H. Weber) (Birhhäuser, Basel, 1980) 103.

[16] G. Iooss and W.F. Langford, Ann. N.Y. Acad. Sci. 357 (1980) 489.

# FIGURE CAPTIONS

**Figure 1: Symmetry Breaking from Basic State.** (a) Basic state, radial forcing, O(2) symmetry; (b) Mode $m = 3$, three pairs of convection rolls; (c) Mode $m = 4$, four pairs of convection rolls.

**Figure 2: Neutral Stability Curves.** Codimension-one curves intersect at codimension-two points in the (Ra,$\eta$) plane, where Ra is Rayleigh number and $\eta$ is radius ratio.

**Figure 3: Sample Bifurcation Diagrams for Large Mode Number.** Both primary bifurcations are assumed supercritical. The branches are labeled by their isotropy groups. Stable branches are indicated by heavy lines.

**Figure 4: Sample Bifurcation Diagrams for Small Mode Number.** One primary branch is assumed subcritical and the other one supercritical. The labels SW and RW denote standing waves and rotating waves respectively and the other branches are labeled by their isotropy groups.

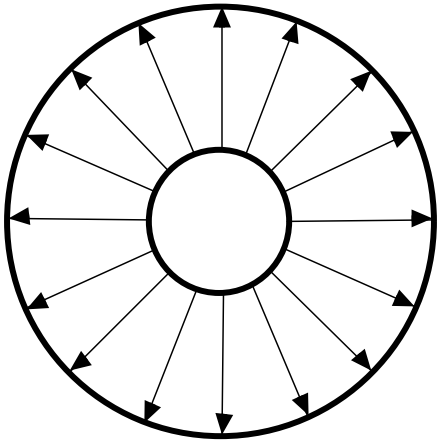 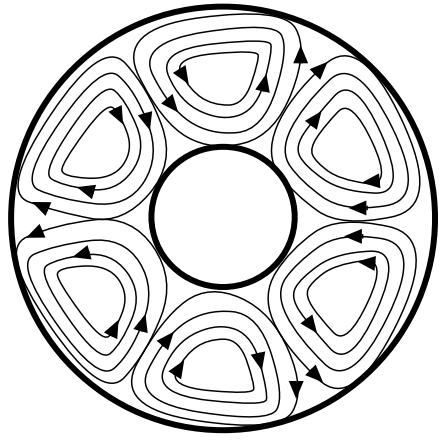 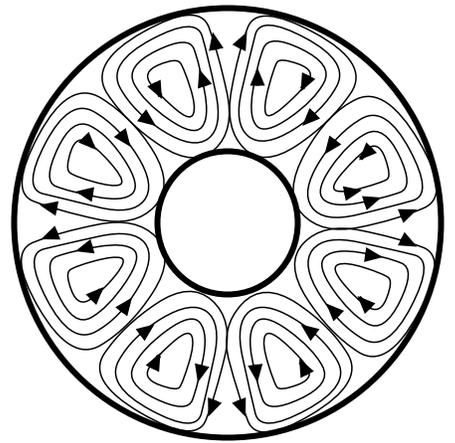

(a)                  (b)                  (c)

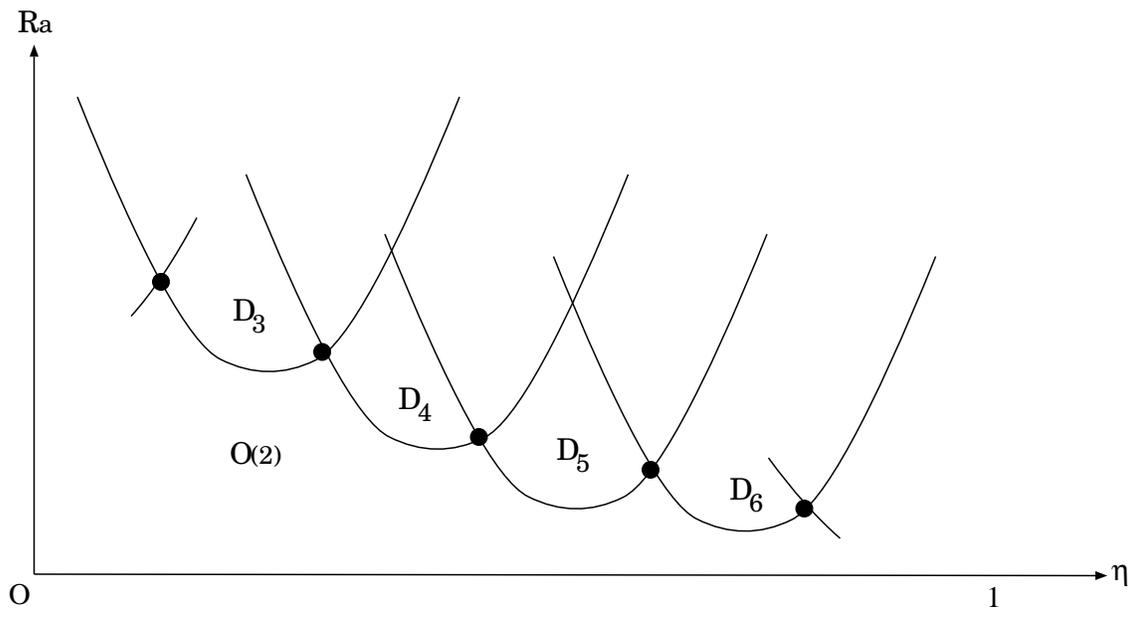

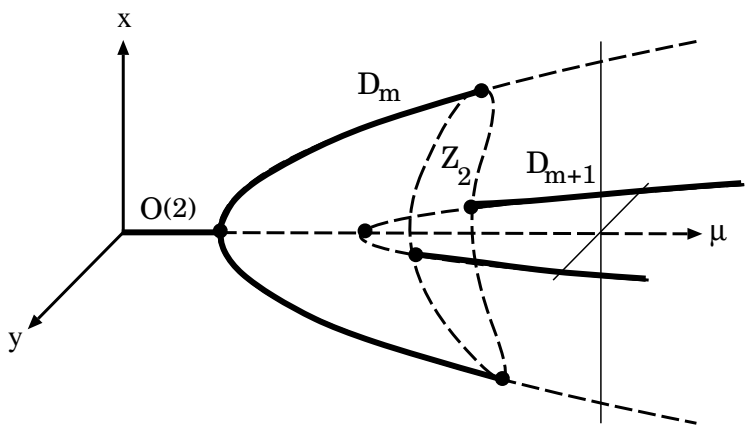
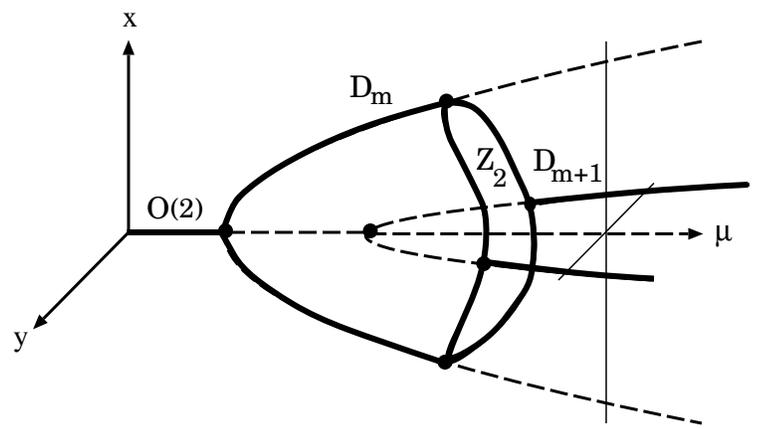
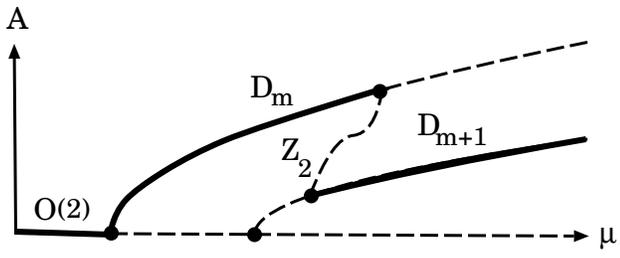
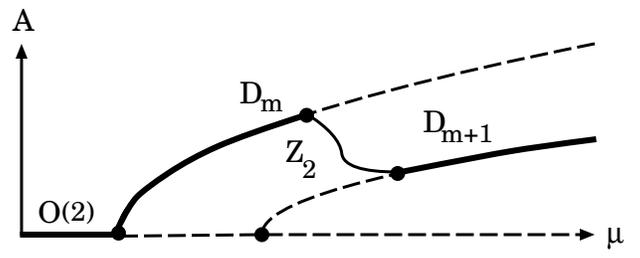

(a)                        (b)

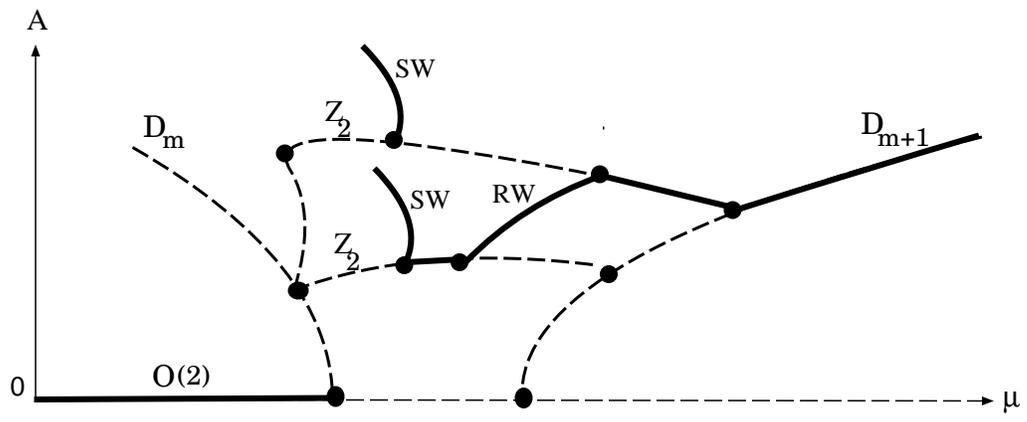